\documentclass[11pt]{article}
\usepackage{amsmath}
\usepackage{amsthm}
\usepackage{amssymb}
\usepackage{amscd}
 \makeatletter

 \@addtoreset{equation}{section}
  \makeatother

\theoremstyle{plain}
\newtheorem{thm}{Theorem}[section]
\newtheorem{prop}[thm]{Proposition}

\newtheorem{lem}[thm]{Lemma}
\theoremstyle{definition}

\theoremstyle{remark}
\newtheorem{rem}{Remark}[section]

\begin{document}
\title{  Chern-Simons pre-quantizations over four-manifolds }
\author{ Tosiaki Kori\thanks
{Research supported by
Promotion for Sciences of the Ministry of
 Education and Science in Japan ( no. 16540084 )}\\
Department of Mathematics\\
School of Science and Engineering\\
 Waseda University \\
3-4-1 Okubo, Shinjuku-ku
 Tokyo, Japan.\\ e-mail: kori@waseda.jp
}
\date{ }
\maketitle

\begin{abstract} 
   We endow the space  of connections  on an \(SU(n)\)-principal bundle over a four-manifold with a pre-symplectic structure and define a Hamiltonian action on it of the group of gauge transformations that are trivial on the boundary.   Then we consider the trivial \(SU(n)\)-principal bundle for \(n\geq 3\) over the four-manifold that is a submanifold of a null-cobordant four-manifold, and  we construct on the  moduli space of connections, as  well as on that of flat connections, a hermitian line bundle with connection whose curvature is given by the pre-symplectic form.    This is the Chern-Simons pre-quantization of moduli spaces.   The group of gauge transformations on the boundary of the four-manifold acts on the moduli space of flat connections by an infinitesimally symplectic way.      When the four-manifold is a 4-dimensional disc we show that this action is lifted to the pre-quantization by its Lie group extension.   The geometric description of the latter is related to the 4-dimensional Wess-Zumino-Witten model.
\end{abstract}


Subj. Class: Global analysis, Quantum field theory.

   Moduli space of flat connections.  Chern-Simons functionals,   Wess-Zumino-Witten actions.

\section{Introduction}

The moduli space of flat connections on a Riemann surface has been investigated by many physicists and mathematicians when its relation to the Chern-Simons theory had been discovered.   It is a compact, finite-dimensional symplectic space.      Ramadas, Singer and Weitsman described the Chern-Simons pre-quantization of this moduli space ~\cite{A, AB,  RSD} .      On the other hand, Donaldson proved that if the surface has a boundary then the moduli space of flat connections is a smooth infinite-dimensional symplectic manifold, and it has a  Hamiltonian action of the group of gauge transformations on the boundary~\cite{D}.    Meinrenken and Woodward~\cite{Me} gave its Chern-Simons pre-quantization.   
In general the moduli space of the \(Lie(G)\)-valued flat connections on a  connected closed manifold \(M\) corresponds bijectively to the conjugate classes of the \(G\)-representations of the fundamental group \(\pi_1(M)\).   

In this paper we study the  Chern-Simons pre-quantization of the moduli space of flat connections on a four-manifold \(M\) generally with non-empty boundary.     Let \({\cal A}(M)\) be the space of irreducible connections on the trivial \(SU(n)\)-principal bundle on \(M\)  and let \({\cal G}_0(M)\) be the group of gauge transformations on \(M\) that are {\it identity on the boundary} \(\,\partial M\).    We shall prove in section 2 that \({\cal A}(M)\) carries a pre-symplectic structure.    The pre-symplectic form is given by 
\begin{equation}
\omega_A(a,b)=\frac{1}{8\pi^3}\int_MTr[\,(ab-ba)F_A\,] -\frac{1}{24\pi^3}\int_{\partial M}Tr[\,(ab-ba)\,A\,],
\end{equation}
for \(a,\,b\in T_A{\cal A}\).     The action of \({\cal G}_0(M)\) becomes a Hamiltonian action with the moment map given by the square of curvature \(F_A^2\).     We can not expect that the pre-symplectic reduction \(\{F_A^2=0\}/{\cal G}_0(M)\) exists as a smooth manifold but the \(0\)-level set of the moment map contains  as its subspace the set of flat connections \({\cal A}^{\flat}(M)\). 
The moduli space of flat connections  \({\cal M}^{\flat}( M)={\cal A}^{\flat}(M)/{\cal G}_0(M)\) is a smooth manifold if  \(\partial M\neq\emptyset\,\).   In this case  
\({\cal M}^{\flat}( M)\) is endowed with a pre-symplectic form 
\begin{equation}
\omega^{\flat}_{[A]}(a,b)= -\frac{1}{24\pi^3}\int_{\partial M}Tr[\,(ab-ba)\,A\,] \,.\end{equation}

   The concept of {\it pre-quantization } of a symplectic manifold is generalized by Guillemin et al. \cite{GGK} to a manifold \(X\) equipped with a closed 2-form ( pre-symplectic form) \(\sigma\).   A pre-quantization of \((X,\sigma)\) is a hermitian line bundle \((L,\,<\,,\,>)\) equipped with a Hermitian connection \(\nabla\) whose curvature is \(\sigma\).   We treat the case where \(X\) is a moduli space of connections on a manifold and \(\sigma\) is a closed 2-form on it.    We shall call this pre-quantization a {\it Chern-Simons pre-quantization} if the transition functions of the quantization \((L,\,<\,,\,>)\) come from Chern-Simons form.   The pre-quantization of the moduli space of flat connections over a Riemann surface was first obtained by ~\cite{AB} and developped by many authors to various directions.  In this case the moduli space is endowed with a symplectic form.      It was pointed out in~\cite{RSD}  that, when the base manifold \(M\)  is a Riemann surface the variation of Chern-Simons 3-form under gauge transformations provides also the transition functions of a hermitian line bundle with connection over the total moduli space and that the pre-quantization of the moduli space of flat connections is obtained by restricting it to flat connections.   In section 3 we extend these results to \(SU(n)\)-bundles with \(n\geq 3\) on four-manifolds.   We construct a hermitian line bundle \(\pi:\,{\cal L}(M)\longrightarrow {\cal B}(M)={\cal A}(M)/{\cal G}_0(M)\) whose transition functions come from 5-dimensional Chern-Simons form.    For that 
we assume after section 3 that the four-manifold \(M\) is a submanifold of a closed four-manifold \(\widehat M\) which is null-cobordant, that is, \(\widehat M\) is the boundary of an oriented 5-dimensional manifold \(N\).    
As a particular case when \(M\) is without boundary; \(M=\widehat M\), we have the line bundle \({\cal L}(\widehat M)\) over \({\cal B}(\widehat M)\).     The 5-dimensional Chern-Simons functional  \({\rm CS}({\bf A})\) of  \({\bf A}\in {\cal A}(N)\) gives a section of the pullback line bundle \(r^{\ast}{\cal L}(\widehat M)\) by the boundary restriction map \(r:{\cal B}(N)\longrightarrow {\cal B}(\widehat M)\).   The gradient vector field of this section is
\[\nabla \,s({\bf A})\,=\frac{i}{8\pi^3}\ast{\bf F}^2_{\bf A},\]
where \(\nabla\) is the covariant differentiation associated to the pullback of the connection on \({\cal L}(\widehat M)\).  
  In section 4 we investigate the pre-quantization of the moduli space \(({\cal M}^{\flat}(M),\omega^{\flat})\).   
The pre-quantization \({\cal L}^{\flat}(M)\) is defined as the  restriction of  \({\cal L}(M)\) to \({\cal M}^{\flat}( M)\subset {\cal B}(M)\).      For a proper submanifold \(M\subset \widehat M\) it is a line bundle with connection whose curvature is  \(-i\,\omega^{\flat}\).   On the other hand the line bundle  \({\cal L}^{\flat}(\widehat M)\) admits a flat connection.     
When  \(\widehat M\) is simply connected \( {\cal M}^{\flat}(\widehat M)\) becomes one-point and \({\cal L}^{\flat}(\widehat M)\simeq {\rm C}\).     

  Other than \({\cal G}_0(M)\);  the group of gauge transformations that are identity on the boundary, we have 
the group of total gauge transformations  \({\cal G}(M)\).   It acts on 
\({\cal A}^{\flat}(M)\) as infinitesimal symplectic automorphisms.   So is the action of \({\cal G}(M)/{\cal G}_0(M)\) on \({\cal M}^{\flat}(M)\).     When \(M\) is the 4-dimensional disc \(D\) with boundary \(S^3\),  \( {\cal G}(D)/{\cal G}_0(D)\) becomes the set of pointed mappings from \(S^3\) to \(G\) that admit extensions to \(D\), which we denote by \(\Omega^3_0G\).      The  infinitesimal symplectic action of \(\Omega^3_0G\) on \({\cal M}^{\flat}(D)\) does not lift to the pre-quantum line bundle \({\cal L}^{\flat}(D)\).   
    Let \(\widehat {\Omega G} \) be the Mickelsson's  {\it abelian extension} of   \(\Omega^3_0G\) by the group \(Map({\cal A}_3,\,U(1))\), \cite{Mi}.     Then  \(\widehat {\Omega G} \) acts on  \({\cal L}^{\flat}(D)\)   and the action 
is equivariant with respect to the action of \(\Omega^3_0G\) on the base space 
 \({\cal M}^{\flat}(D)\).   The reduction of  \({\cal L}^{\flat}(D)\) by the action  of  \(\widehat {\Omega G} \) becomes \({\rm C}\,\).    
These are the analogy of the central extension of the loop group acting on the quantization of the moduli space of flat connections over a surface~\cite{Me}.    The associated line bundle to the abelian extension 
\(\widehat{\Omega G}\) is isomorphic to the restriction of \({\cal L}^{\flat}(D)\) to the space of pure gauges \( \{g^{-1}dg \in {\cal M}^{\flat}(D):\,g\in {\cal G}(D)\}\), which we denoted in \cite{K} as \(WZ(S^3)\).   We extend the definition \(WZ(\Gamma)\) to the disjoint union \(\Gamma=\bigoplus S^3\).   Let \(\Sigma\) be a four-manifold that cobords \(\Gamma\).    We can assign to \(\Sigma\) a trivialization \(WZ(\Sigma)\) of the pullback bundle \(r^{\ast}WZ(\Gamma)\), \(r\) being the boundary restriction.   Then \(WZ\) gives a functor that assigns to the category of three manifolds \(\Gamma\) and their cobordants \(\Sigma\) the category of the line bundles over the space of mappings \(Map_0(\Gamma,G)\).    \(WZ\) satisfies the axioms of Atiyah-Segal's TQFT.   By axiomatizing the multiplication properties satisfied by \(WZ(\Sigma)\) we introduced in \cite{K} 4-dimensional Wess-Zumino-Witten theory.

\section{ Preliminaries}

\subsection{Differential calculations on the space of connections}  

   Let \(M\) be an oriented Riemannian four-manifold with boundary \(\partial M\).   Let \(G=SU(n)\).   The inner product on \(G\) is given by 
\(
<\xi,\eta>=
-Tr(\xi\eta)\) for \( \xi,\,\eta\in Lie(G)=su(n)\).    
With this inner product the dual of \(Lie (G)\) is identified with \(Lie (G)\) itself.      
Let \(P= M\times G\) be the trivial principal \(G\)-bundle over \(M\).   
We define \(L^2_{s-1}\) connections on \(P\)
by \(L^2_{s-1}\) connection matrices.   We write \({\cal A}={\cal A}(M)\) for the space of {\it irreducible} \(L^2_{s-1}\) connections.   
The tangent space of \({\cal A}\) at \(A\in{\cal A}\) is \(T_A{\cal A}=\Omega^1_{s-1}(M, Lie\,G)\).    
The curvature of \(A\in {\cal A}\) is 
\(F_A=dA+\frac12[A\wedge A]\in \Omega^2_{s-2}(M,  Lie\,G)\).   
\(\Omega^3_{s-1}(M,  Lie\,G)\), being the formal dual of \(\Omega^1_{s-1}(M,  Lie\,G)=T_A{\cal A}\,\), is identified with the space of 1-forms on \({\cal A}\).

Here are some differential calculations on \({\cal A}\) that we shall cite from~\cite{BN, DK,S}.   
The derivation of a smooth function \(H=H(A)\) on \({\cal A}\) is defined by the functional derivation of \(A\):
\begin{equation}
(\partial_A H)a=\lim_{t\longrightarrow 0}\frac{H(A+ta)-H(A)}{t},\quad \mbox{ for \(a\in T_A{\cal A}\)}.
\end{equation}
For example, 
\((\partial_A A)a=a\), and \((\partial_A F_A)a=d_Aa\).
The second follows from the formula
\(F_{A+a}=F_A+d_Aa+a\wedge a\).   
The derivation of a vector field on \({\cal A}\)  is  defined as that of a smooth function of \(A\in {\cal A}\) valued in \(\Omega^1_{s-1}(M,Lie\,G)\), similarly the derivation of  a 1-form on \({\cal A}\)  is done in \(\Omega^3_{s-1}(M,Lie\,G)\).   
We have, for a vector field \({\bf b}\) and a 1-form \(\beta\), 
\begin{equation}
(\partial_A<\beta, {\bf b}>)a=<\beta,\,(\partial_A{\bf b})a>+<(\partial_A\beta)a,{\bf b}>.  \label{Leibnitz}
\end{equation}
The Lie bracket for vector fields on \({\cal A}\) is seen to have the expression
\begin{equation}
[{\bf a},\,{\bf b}]=(\partial_A{\bf b}){\bf a}-(\partial_A{\bf a}){\bf b}. \label{Liebra}
\end{equation}
Let \(\widetilde d\)  be the exterior derivative on \({\cal A}\).   For a function \(H=H(A)\) on \({\cal A}\), it is defined by  \((\widetilde dH)_Aa=(\partial_AH)a\).   
 (\ref{Leibnitz} ) and (\ref{Liebra}) yield the following formula for the exterior derivative of a 1-form \(\theta\) on \({\cal A}\) :
\begin{equation}
(\widetilde d\theta)_A({\bf a},{\bf b})= <(\partial_A\theta){\bf a},  {\bf b}>-<(\partial_A\theta){\bf b},{\bf a}>.  \label{extder1}
\end{equation}
For a 2-form \(\varphi\), 
\begin{equation}
(\widetilde d\,\varphi)_A({\bf a},{\bf b},{\bf c})=(\partial_A\varphi({\bf b},{\bf c})){\bf a}+(\partial_A\varphi({\bf c},{\bf a})){\bf b}+(\partial_A\varphi ({\bf a},{\bf b})){\bf c}.   \label{extder2}  \end{equation}

\subsection{ Moduli space of connections }

We denote \({\cal G}={\cal G}(M)\) the group of \(L^2_s\)-gauge transformations based at a point \(p_0\in M\):
\begin{equation}
{\cal G}(M)=\, \{g\in \Omega^0_s(M, G)\,;\,g(p_0)=1\,\}.\end{equation}   
\({\cal G}\) acts on \({\cal A}\) by 
\(
g\cdot A=g^{-1}dg+g^{-1}Ag=A+g^{-1}d_Ag
\).    The action being free,  the orbit space  \({\cal A}/{\cal G}\) is a smooth infinite dimensional manifold.    By Sobolev lemma one sees that \({\cal G}\) is a Banach Lie group  and its action is a smooth map of Banach manifolds.   We have  \(Lie({\cal G})=\Omega^0_s(M, Lie\,G)\).

When \(M\) has the boundary we assume \(p_0\in \partial M\), and the group of \(L^2_{s-\frac12}\)-gauge transformations on the boundary \(\partial M\) is denoted by \({\cal G}(\partial M)\).   We have the restriction map to the boundary:
\(r:\,{\cal G}(M)\longrightarrow {\cal G}(\partial M)\).   
Let \({\cal G}_0={\cal G}_0(M)\) be the kernel of the restriction map.   
It is the group of gauge transformations that are identity on the boundary.   When \(\partial M\neq\emptyset\) the moduli space 
 \({\cal A}/{\cal G}_0\) is a smooth infinite dimensional manifold.   
  
 The derivative of the action of \(\,{\cal G}\) at \(A\in{\cal A}\) is 
\begin{equation}
d_A=d+[A\wedge\quad]:\, \Omega^0_s(M, Lie\,G)\longrightarrow \Omega^1_{s-1}(M,Lie\,G).
\end{equation}
The fundamental vector field on \({\cal A}\) corresponding to  
\(\xi \in Lie({\cal G})\) is given by 
\(\xi_{{\cal A}}(A)=\frac{d}{dt}\lvert_{t=0}(\exp \,t\xi)\cdot A=d_A\xi\),   
and the tangent space to the orbit at \(A\in{\cal A}\) is 
\begin{equation}
T_A({\cal G}\cdot A)=d_A\,Lie({\cal G}).
\end{equation}
By Stokes'  formula we have  the orthogonal decomposition
\begin{eqnarray}
T_A{\cal A}&=&\,d_A\, Lie({\cal G})\,\oplus \,H_A,\\[0.2cm]
H_A&=&\{a\in\Omega^1_{s-1}(M, Lie\,G);\, \,d_A^{\ast}a=0,\,\mbox{and } \,\ast a\vert\partial M=0\}.\nonumber
\end{eqnarray}
As for the orbit of   \({\cal G}_0\), we have the orthogonal decomposition:
\begin{eqnarray}
T_A{\cal A}&=&\,d_A\, Lie({\cal G}_0)\,\,\oplus \,H^0_A,\\[0.2cm]
H^0_A&=&\{a\in\Omega^1_{s-1}(M, Lie\,G);\, d_A^{\ast}a=0\}.\nonumber
\label{orthodecomposition}
\end{eqnarray}

We have two moduli spaces of irreducible connections;
\begin{equation}
 {\cal B}(M)={\cal A}/{\cal G}_0 ,\qquad {\cal C}(M)={\cal A}/{\cal G}.
\end{equation}
\({\cal B}(M)\) is a \({\cal G}/{\cal G}_0\)-principal bundle over \({\cal C}(M)\).    \({\cal C}(M)\) coincides with \({\cal B}(M)\) if  \(M\) has no boundary.     \({\cal C}(M)\) is finite dimensional but in general \({\cal B}(M)\) is infinite dimensional, in fact it contains the orbit of \({\cal G}(\partial M)\).
   \(  {\cal B}(M)\) is a smooth manifold modelled locally on the balls in the Hilbert space \(\ker d^{\ast}_A\)  in  \(\Omega^2_{s-1}(M, G)\).       \({\cal C}(M)\)  is a smooth manifold modelled locally on the balls in the Hilbert space \( \ker d^{\ast}_A \cap \ker(\ast\vert\partial M)\) in  \(\Omega^2_{s-1}(M, G)\).     The reader can find the precise and technical description of these facts in ~\cite{D, DK}. 
   
   The space of flat connections is 
\begin{equation}
{\cal A}^{\flat}(M)=\{A\in{\cal A}(M); F_A=0\},
\end{equation}
which we shall often abbreviate to \({\cal A}^{\flat}\).   
The tangent space of \({\cal A}^{\flat}\) is given by
\begin{equation}
T_A{\cal A}^{\flat}=\{a\in \Omega^1_{s-1}(M, Lie\,G); \, d_Aa=0\}.
\end{equation}
Any vector  tangent to the \({\cal G}\)-orbit through \(A\in{\cal A}^{\flat}\) is in \(T_A{\cal A}^{\flat}\,\).   We have the orthogonal decomposition 
\begin{equation}
T_A{\cal A}^{\flat}=d_ALie({\cal G}_0)\oplus H^{0\,\flat}_A,\label{decomp}
\end{equation}
with
\[ H^{0\,\flat}_A=\{\Omega^1_{s-1}(M, Lie\,G );\, \,d_Aa=0,\,d^{\ast}_Aa=0\}.\]
The moduli space of flat connections is by definition 
\begin{equation}
{\cal M}^{\flat}={\cal A}^{\flat}/{\cal G}_0 .\label{flatmoduli1}
\end{equation}
  We know that  \({\cal M}^{\flat}\) is a smooth manifold~\cite{DK}.     
On the other hand 
we have the following orthogonal decomposition by the action of the total gauge transformation group \({\cal G}\):
\begin{equation}
T_A{\cal A}^{\flat}=d_ALie({\cal G})\oplus H^{\flat}_A,
\end{equation}
with
\[ H_A^{\flat}=\{a\in\Omega^1_{s-1}(M,  Lie\,G);\, \,d_Aa=d_A^{\ast}a=0,\,\mbox{and } \ast a\vert\partial M=0\}.\]
The moduli space of flat connections modulo the group \({\cal G}\) is 
\begin{equation}
{\cal N}^{\flat}={\cal A}^{\flat}/{\cal G}.\label{flatmoduli2}
\end{equation}
The dimension of \({\cal N}^{\flat}\) is finite.

\section{Pre-symplectic structure on \({\cal A}\)}

\subsection{Pre-symplectic structure on \({\cal A}\) and the action of  \({\cal G}_0(M)\)}

 In the sequel we shall suppress the Sobolev indices.     So \({\cal A}\) is always the space of irreducible \(L^2_{s-1}\)-connections and \({\cal G}\) is the group of based \(L^2_s\)-gauge transformations.   
 
For each \(A\in {\cal A}\) we define a sqew-symmetric bilinear form 
on \(T_A{\cal A}\) by:
\begin{eqnarray}
\omega_A(a,b)&=&\omega^0_A(a,b)+\omega^{\prime}_A(a,b),\\ [0.2cm]
\omega_A^0(a,b)&=&\frac{1}{8\pi^3}\int_M\,Tr[\,(a\wedge b-b\wedge a)\wedge F_A\,],\\[0.2cm]
\omega^{\prime}_A(a,b)&=&-\frac{1}{24\pi^3}\int_{\partial M}\,Tr[\,(a\wedge b-b\wedge a)\wedge A\,],
\end{eqnarray}
for \(a,b\in T_A{\cal A}\).   

\begin{thm}
 \(\omega\) is a \({\cal G}_0\)-invariant closed 2-form on \({\cal A}\).
\end{thm}
{\it Proof.}~~
\({\cal G}_0\) invariance of \(\omega\) is evident since any \(g\in{\cal G}_0\) is identity on \(\partial M\).   
In the following we shall abbreviate \(ab\) for the exterior product \(a\wedge b\).   
By the formula (\ref{extder2}), we have, for \(a,b,c\in T_A{\cal A}\), 
\begin{equation*}
(\widetilde d\,\omega^0)_A(a,b,c)=
(\partial_A \omega^0(a,b))c+ (\partial_A \omega^0(b,c))a+ ( \partial_A \omega^0(c,a))b\,,.
\end{equation*}
Since
\[(\partial_A \omega^0(a,b))c=\frac{1}{8\pi^3}\int_M\,Tr[(ab-ba)d_Ac],
\]
\begin{eqnarray*}
(\widetilde d \omega^0)_A(a,b,c)&=&\frac{1}{8\pi^3}\int_M\,Tr\left [\,(ab-ba)d_Ac+
(bc-cb)d_Aa+(ca-ac)d_Ab\right]\\[0.2cm]
&=&\frac{1}{8\pi^3}\int_Md\,Tr[\,(ab-ba) c\,]=\frac{1}{8\pi^3}
\int_{\partial M}Tr[(ab-ba)c].\end{eqnarray*}
On the other hand we have
\[(\widetilde d \omega^{\prime})_A(a,b,c)=3(\partial_A \omega^{\prime}(a,b))c=-\frac{1}{8\pi^3}\int_{\partial M}Tr[(ab-ba)c].\]
Therefore
\(\widetilde d\, \omega=0\).
\hfill \qed

\begin{rem}   
For \(M\) without boundary the pre-symplectic form \(\omega^0\) was introduced by Bao and Nair~\cite{BN}.    \end{rem}

Moreover we find that \(\omega\) is an exact form.   
\begin{lem}\label{exact}
Let \(\theta\)  be a 1-form on \({\cal A}\) defined by
\begin{equation}
\theta_A(a)=-\frac{i}{24\pi^3}\int_M\,Tr[\, (A F+F A-\frac12 A^3)\, a \,],\quad a\in T_A{\cal A}.\label{theta}
\end{equation}
Then
\begin{equation}
\widetilde d\,\theta\,=\,-i\omega.
\end{equation}
\end{lem}
{\it Proof.}~~
For \(a,b \in T_A{\cal A}\), 
\begin{eqnarray*}
(\widetilde d\theta)_A(a,b)&=&\langle (\partial_A\theta)a,b\rangle - \langle (\partial_A\theta)b, a\rangle \\[0.2cm]
&=&  \frac{-i}{24\pi^3}\int_M\,Tr[\,2(ab-ba)F-(ab-ba)A^2\\[0.2cm]
&& \qquad\qquad - \,(bd_Aa+ d_Aab-d_Aba-ad_Ab)A\,] .
\end{eqnarray*}
But since
\[d\,Tr[(ab-ba)A]=Tr[(b\,d_Aa+ d_Aa\,b-d_Ab\,a-a\,d_Ab)A]+Tr[(ab-ba)(F+A^2)],\]
we have 
\begin{eqnarray*}
(\widetilde d\,\theta)_A(a,b)&=& - \frac{i}{8\pi^3}\int_MTr[(ab-ba)F ] +\frac{i}{24\pi^3}\int_{\partial M}Tr[(ab-ba)A]\nonumber\\
&=&
-i\, \omega_A(a,b).
\end{eqnarray*}
\hfill\qed

We note that \(\theta\) is not \({\cal G}_0\) invariant.

\subsection{Hamiltonian \({\cal G}_0-\)manifold \(({\cal A},\omega,\Phi)\)}

We shall explain briefly the generalization by Guillemin et al. \cite{GGK} of the notion of Poisson algebra of a symplectic manifold to manifolds \((X,\sigma)\) equipped with closed 2-forms. 

Let \(X\) be a manifold with a smooth action on it of a Lie group \(G\).   Let  \(\sigma\) be a \(G\)-invariant closed 2-form on \(X\).  For any \(\xi\in Lie\,G\) there corresponds a vector field \(\xi_X\) generating the action on \(X\) of 1-parameter subgroup \(\{\exp\,t\xi;\,t\in {\rm R}\}\subset G\).   
A moment map \(\Phi:X\longrightarrow (Lie\,G)^{\ast}\) is a map that is equivariant with respect to the \(G\)-action on \(X\) and the coadjoint action on \((Lie\,G)^{\ast}\) such that \(\Phi^{\xi}=<
\Phi,\xi>\) satisfies 
\[d\Phi^{\xi}=i(\xi_X)\sigma,\quad \forall\xi\in Lie\,G.\]
We call \(\left(X,\sigma,\Phi\right)\) a {\it Hamiltonian \(G\)-manifold}. 
Put
\begin{equation}
{\cal F}(X , \sigma)=\{(f,v)\in C^{\infty}(X)\times Vect(X):\quad df=i_v\sigma\,\}.
\end{equation}
We endow \({\cal F}(X ,\sigma)\)  the multiplication defined by
\begin{equation}
(f,v)\cdot (g,u)=(fg,fu+gv),
\end{equation}
and the bracket defined by
\begin{equation}
\left[(f,v),(g,u)\right]=\left(L_uf,-[u,v]\right),
\end{equation}
where
\[L_uf=-L_vg=\sigma(u,v).\]
\begin{prop}[\cite{GGK}]
\({\cal F}(X ,\sigma)\) is a Poisson algebra, that is, it is a commutative algebra, it is a Lie algebra and the Lie bracket is a derivation with respect to the multiplication.
\end{prop}
\begin{prop}[\cite{GGK}]
For any Hamiltonian \(G\)-manifold \((X,\sigma,\Phi)\), the map
\begin{equation}
Lie\, G\ni\xi \longrightarrow (\,\Phi^{\xi},\,\xi_X\,)\in {\cal F}(X,\sigma)
\end{equation}
gives a homomorphism of Lie algebras.
\end{prop}
The proof of these propositions are standard.

  The above mentioned notion of Poisson algebra fits well to our case if we consider our Hamiltonian \({\cal G}_0\)-manifold \(({\cal A},\omega,\Phi)\) that we shall explain in the sequel.   We have already shown that \(({\cal A}, \omega)\) is a manifold equipped with a \({\cal G}_0\)-invariant closed 2-form.    The moment map will be given by the following theorem and we have the Hamiltonian \({\cal G}_0\)-manifold \(({\cal A},\omega,\Phi)\).

\begin{thm}~
The action of \({\cal G}_0\) on \({\cal A}\) is a Hamiltonian action and the corresponding moment map is given by 
\begin{eqnarray}
\Phi&:&{\cal A}\longrightarrow (Lie\, {\cal G}_0)^{\ast}=\Omega^4(M, Lie\,G):\,\quad A\longrightarrow F_A^2\,.\nonumber\\
\nonumber\\
\langle \Phi(A), \xi\rangle&=&\Phi^{\xi}(A)\,=\,\frac{1}{8\pi^3}\int_M\,Tr(F_A^2\xi),\, \mbox{ for }
\xi\in Lie\,{\cal G}_0\,.
\end{eqnarray}
\end{thm}
{\it Proof.}~~
The equivariance of \(\Phi:{\cal A}\longrightarrow (Lie\, {\cal G}_0)^{\ast}\) with respect to the \({\cal G}_0\)-action on \({\cal A}\) and the coadjoint action on \((Lie\, {\cal G}_0)^{\ast}\) is evident.     We have 
\[(\partial \Phi^{\xi})_Aa=\frac{1}{8\pi^3}\int_M\,Tr[(d_Aa\wedge F_A+F_A\wedge d_Aa)\xi],\]
and
\begin{equation*}
\begin{split}
(\partial \Phi^{\xi})_Aa&-\omega^0_A(a,d_A\xi)=
\frac{1}{8\pi^3}
\int_M\,Tr[(d_Aa\wedge F_A+F_A\wedge d_Aa)\xi-(F_Aa+aF_A)d_A\xi]\\[0.2cm]
&=\frac{1}{8\pi^3}\int_MdTr[(aF_A+F_Aa)\xi]=\frac{1}{8\pi^3}\int_{\partial M}Tr[(aF_A+F_Aa)\xi].
\end{split}
\end{equation*}
It is equal to  \(0\) since \(\xi=0\) on \(\partial M\).
   On the other hand it holds that
 \begin{eqnarray*}
dTr[(Aa-aA)\xi]&=&Tr[(F_Aa+aF_A-Ad_Aa-d_AaA+A^2a+aA^2)\xi]\\[0.2cm]
&&\qquad \qquad \qquad +\, Tr[ (ad_A\xi-d_A\xi a)A]
,\end{eqnarray*}
the first term of which vanishes on \(\partial M\).
 Hence  
\begin{equation}\label{bomega}
\omega^{\prime}(a,d_A\xi)=-\frac{1}{24\pi^3}\int_{\partial M}
dTr[(Aa-aA)\xi]=0.\end{equation}
Therefore
\[ (\partial \Phi^{\xi})_Aa=\omega_A(a,d_A\xi).\]
\hfill\qed 
\begin{prop}
We have, for  \(\xi, \,\eta\in Lie\, {\cal G}_0\,\), 
\begin{equation}
\Phi^{[\xi,\eta]}(A)=\omega_A(\,d_A\xi,\,d_A\eta\,).
\end{equation}
   \end{prop}
{\it Proof} ~~ 
For \(\xi, \,\eta\in Lie\, {\cal G}_0\), we have
\[
d\,Tr[\,(\xi\,d_A\eta-d_A\xi\,\eta ) F_A\,] = Tr[\,2d_A\xi\,d_A\eta\, F_A - \xi\,\eta F_A^2\,+\,\xi F_A\,\eta F_A\,].
\]
This equation and the same one with \(\xi\) and \(\eta\) reversed yield
\[ 2\,Tr\left[\,(d_A\xi\,d_A\eta-d_A\eta\,d_A\xi)F_A\,-\, [\xi,\,\eta]\,F_A^2\,\right]=\,d\,Tr[\,\cdots\,]\,.\]
Integration over \(M\) yields
\[\int_M\,\,Tr \left(\,F_A^2\,[\xi,\,\eta]\,\right)=\int_M\,Tr\left[\,(d_A\xi\,d_A\eta-d_A\eta\,d_A\xi)F_A\,\right].\]
Since \(\omega^{\prime}_A(d_A\xi,\,d_A\eta)=0\) from (\ref{bomega}) we have the assertion.

\hfill\qed 

\begin{prop}
The map \(\Phi:\, Lie\,{\cal G}_0\ni\xi \longrightarrow (\Phi^{\xi},\,\xi_{{\cal A}})\in {\cal F}({\cal A},\omega)\)  gives a homomorphism of Lie algebras,  
where \(\xi_{{\cal A}}(A)=d_A\xi\).
\end{prop}

\subsection{Reduction to the moduli space of flat connections}

We have 
\(\Phi^{-1}(0)=\{A\in {\cal A};\, F_A^2=0\}\) .      
Since \(g^{\ast}\omega=\omega\) for \(g\in {\cal G}_0\,\) and,   \(i_{d_A\xi}\omega=0\) on \(\Phi^{-1}(0)\) for \(\xi\in Lie\,{\cal G}_0\,\),  
the 2-form \(\omega\) descends to \(\Phi^{-1}(0)/{\cal G}_0\) and gives a closed 2-form on it.   This is the  pre-symplectic reduction if \(\Phi^{-1}(0)\) is a manifold and the action of \({\cal G}_0\) on it is locally free.         The 
\({\cal G}_0\)- action is indeed free but \(0\) is not a regular value of the moment map \(\Phi\), so we can not expect that \(\Phi^{-1}(0)\) is a smooth manifold.   But  it contains the space of flat connections  \({\cal A}^{\flat}\), which is known to be a smooth manifold if \(\partial M\neq\emptyset\) .     
\begin{thm}\label{bsymplectic}~\
Suppose \(\partial M\neq\emptyset\).   Then 
\({\cal M}^{\flat}\) is a smooth manifold endowed with a pre-symplectic structure \(\omega^{\flat}\).   
The pre-symplectic form on \({\cal M}^{\flat}\) is given by 
\begin{equation}
\omega^{\flat}_{[A]}(a, b)=\omega^{\prime}_A(a,b)
\end{equation}
for \([A]\in{\cal M}^{\flat}\) and \(a,\,b \in T_{[A]}{\cal M}^{\flat}\).
\end{thm}
 Here \([A]\in {\cal M}^{\flat}\) denotes the \({\cal G}_0\)-orbit of \(A\in{\cal A}^{\flat}\), and as for a tangent vector \(a\in T_{[A]}{\cal M}^{\flat}\) we take the representative tangent vector to the slice; \(a\in H_A^{0 \,\flat}\), see (\ref{decomp}) .   
\(\omega^{\flat}\) is well defined because we have 
\(g\cdot A=A\) on \(\partial M\) for \(g\in {\cal G}_0\), and \(\omega^{\prime}_A(a,d_A\xi)=0\) for \(\xi\in Lie\,{\cal G}_0\) and \(a\in T_A{\cal A}^{\flat}\), (\ref{bomega}).

\vspace{0.2cm}

{\bf Example 1}.~~
For \(M=S^4\), 
the moduli space of flat connections  
 \[{\cal M}^{\flat}(S^4)={\cal A}^{\flat}(S^4)/{\cal G}(S^4).\]  
is one-point.   \\
In fact, let \(p_0\in S^4\) and let \(A\in {\cal A}^{\flat}(S^4)\).   Let \(T^A_\gamma(x)\) denote the parallel transformation by \(A\) along the curve \(\gamma\) joining \(p_0\) and \(x\).  
We put 
\(t_A(x)=T^A_\gamma(x)1\,\in G\).   It is independent of the choice of curve \(\gamma\) joining \(p_0\) and \(x\).   
 Then \(t_A\in {\cal G}(S^4)\), and by the definition  
\(A=dt_A\cdot t_A^{-1}\).   

In general \({\cal M}^{\flat}(M)\) is one-point for an oriented, connected and simply connected compact four-manifold \(M\) without boundary.\\
{\bf Example 2.}~~
For a disc \(D^4=\{x\in {\rm R}^4;\,\vert x\vert\leq 1\}\) with boundary  \(S^3\), we have
 \begin{equation}{\cal M}^{\flat}(D^4)\simeq \Omega^3_0G.\,\label{OmegaD}
\end{equation}
    Where 
\(\Omega^3G=\{f\in Map(S^3, G); f(p_0)=1\}\), 
and 
\(\Omega^3_0G\) is the connected component of the identity.    
To prove it, first we note that 
\[{\cal G}(D^4)\simeq D^4G=\{f\in Map(D^4,G);\, f(p_0)=1\},\]
and
\[ {\cal G}_0(D^4)\simeq D^4_0G=\{f\in D^4G;\,f\vert S^3=1 \}.\]
Hence \({\cal G}/{\cal G}_0\simeq \Omega^3_0G\).   As before we put, for \(A\in {\cal A}^{\flat}(D^4)\),  \(t_A(x)=T^A_{\gamma}(x)1\,\), \(x\in D^4\).    We have a well defined bijective map from 
\({\cal A}^{\flat}(D^4)\) to \(D^4G\).    In particular, \(t_A=g\) for \(A=dg\,g^{-1}\) with \(g\in {\cal G}(D^4)\) .    
It holds also that \(t_{g\cdot A}(x)=t_A(x)g(x)\) for  \(g\in{\cal G}(D^4)\) .   
Hence we have the isomorphism
\[ {\cal M}^{\flat}(D^4)= {\cal A}^{\flat}/{\cal G}_0\simeq D^4G/D^4_0G\simeq \Omega^3_0G.\]     
 
\subsection{The action of \({\cal G}(M)\)}

By the action of the group of total gauge transformations  \({\cal G}={\cal G}(M)\)  on \({\cal A}^{\flat}\) we have the orbit space \({\cal N}^{\flat}={\cal A}^{\flat}/{\cal G}\).    Then we have a fibration \({\cal M}^{\flat} \longrightarrow {\cal N}^{\flat}\) with the fiber \({\cal G}/{\cal G}_0\), see (\ref{flatmoduli1}) and (\ref{flatmoduli2}).
\begin{prop}
The action of \({\cal G}\) on \({\cal A}^{\flat}\)  is infinitesimally symplectic.   So is the action of 
 \({\cal G}/{\cal G}_0\)  on
 \( ({\cal M}^{\flat}, \omega^{\flat})\) .      
\end{prop}
In fact, we have \(d_A\xi\in T_A{\cal A}^{\flat}\) for \(\xi \in Lie({\cal G})\) and \(A\in {\cal A}^{\flat}\).    Let  \(L_{d_A\xi}\) be the Lie derivative by the fundamental vector field \(d_A\xi\) corresponding to \(\xi\in Lie\,{\cal G}\).   We have, for   
\(a,b\in T_A{\cal A}^{\flat}\),
\begin{equation*}
\begin{split}
(L_{d_A\xi}\omega)_A&(a,b)=(\widetilde d\,i_{d_A\xi}\,\omega)_A(a,b)
=\partial_A\,(i_{d_A\xi}\,\omega_A(b))(a)-\partial_A\,(i_{d_A\xi}\,\omega_A (a))(b)\\[0.2cm]
&= -\frac{1}{24\pi^3}\int_{\partial M}\,Tr[( b\,d_A\xi-d_A\xi\, b)a]+
\frac{1}{24\pi^3}\int_{\partial M}\,Tr[( a\,d_A\xi-d_A\xi \,a)b]\\[0.2cm]
&= -\frac{1}{12\pi^3}\int_{\partial M}\,Tr[(ab-ba)d_A\xi]
= -\frac{1}{12\pi^3}\int_{\partial M}\,d\,Tr[(ab-ba)\xi]\\[0.2cm]
&=0.
\end{split}
\end{equation*} 
Hence \({\cal G}\) acts on \({\cal A}^{\flat}\) by an infinitesimally symplectic way.   
Since \({\cal G}_0\) is a normal subgroup of \({\cal G}\),  the action of 
 \({\cal G}/{\cal G}_0\)  on 
 \( ({\cal M}^{\flat}, \omega^{\flat})\) is  infinitesimally symplectic .      
 \hfill\qed 

\vspace{0.3cm}

{\bf Example 3.}   
The argument in Example 2 shows also that
\[{\cal N}^{\flat}(D^4)=\,\mbox{ one  point} .\]
The same argument by using parallel transformations along the curves in \(S^3\) yields that 
\[\Omega^3G\simeq {\cal A}^{\flat}(S^3).\]    
So we have an injective mapping  \({\cal M}^{\flat}(D^4)\longrightarrow  {\cal A}^{\flat}(S^3)\).   It corresponds to the embedding \(\Omega^3_0G\longrightarrow \Omega^3G\).    

\vspace{1cm}

\section{ Chern-Simons 5-form and the line bundle on \({\cal B}(M)\) }

\subsection{ Chern-Simons 5-form and its variation}

Let  \(G=SU(n)\).   Let \(N\) be an oriented 5- manifold.    As typical examples we are thinking of the 5-dimensional disc \(N=D^5\) and 5-sphere \(S^5\).   Let 
 \(P=N\times G\) be the trivial \(G\)-principal bundle over \(N\).   Let \({\cal A}(N)\) denote the space of connections on \(N\).   We do not suppose that connections on \(N\) are  irreducible.   The group of gauge transformations on \(N\) is denoted by \({\cal G}={\cal G}(N)\).

  Let \(\Omega^q(N)\) be the differential \(q\)-forms on \(N\) and let \(V^q\) be the vector space of polynomials \(\Phi=\Phi(A)\) of \(A\in {\cal A}(N)\)  and its curvature \(F_A\) that take values in \(\Omega^q(N)\).  The curvature \(F_A\) of a connection \(A\) will be often abbreviated to \(F\).   The group of gauge transformations \({\cal G}\) acts on \(V^q\) by \((g\cdot \Phi)(A)=\Phi(g^{-1}\cdot A)\).   We shall investigate the double complex 
\[C^{p,q}=C^p({\cal G}, V^{q+3}),\] 
that is doubly graded by the chain degree \(p\) and the differential form
degree \(q\).   Let \(\,d:C^{p,q}\longrightarrow C^{p,q+1}\) be the exterior differentiation.    The coboundary operator \(\delta: C^{p,q}\longrightarrow C^{p+1,q}\) is given by
\begin{equation*}
\begin{split}
(\delta\, c^p)(g_1,g_2,\cdots,g_{p+1})&=
g_1\cdot c^p(g_2,\cdots,g_{p+1})+(-1)^{p+1}c^p(g_1,g_2,\cdots,g_p)\\[0.2cm]
&+\sum_{k=1}^p(-1)^kc^p(g_1,\cdots,g_{k-1},g_kg_{k+1},g_{k+2},\cdots,g_{p+1}).
\end{split}
\end{equation*}
The following Proposition is a more precise version of the Wess-Zumino's descent equation~\cite{Z}.     The author learned it with its transparent proof from Y. Terashima of Tokyo Institute of Technology.    For \(V^{q+3}\) stated here the calculations were already appeared in~\cite{Mi,Mic}.

\begin{prop}\label{WZdescent}~
There is a sequence of cochains \(c^{p,q}\in C^{p,q}\), \(0\leq p,q\leq 3\),  that satisfies the following relations:
\begin{eqnarray}
dc^{p,3-p}+(-1)^p\delta c^{p-1,3-p+1}&=&0\label{WZdescent1}\\[0.2cm]
dc^{p,2-p}+(-1)^{p}\delta c^{p-1,3-p}&=&-c^{p,3-p}\label{WZdescent2}\\[0.2cm]
c^{0,3}=0,\qquad c^{p,q}&=&0\qquad \mbox{if }\,p+q\neq 2,3 \,.\nonumber
\end{eqnarray}
\end{prop}
   Each term is given by the following forms:
\begin{eqnarray*}
c^{0,2}(A)&=&Tr\,(AF^2-\frac12
A^3F+\frac1{10}A^5),\\[0.2cm]
c^{1,2}(g)&=& \frac1{10}
Tr(dg\cdot g^{-1})^5,\\[0.2cm]
c^{1,1}(g; A)&=& Tr[\,-\frac12V(AF+FA-A^3)+
\frac14(VA)^2+\frac12V^3A\,],\\[0.2cm]
&&\qquad\mbox{where \(V=dg\,g^{-1}\)}, \\[0.2cm]
c^{2,1}(g_1,g_2)&=&c^{1,1}(g_2 ;\,g_1^{-1}dg_1\,),\\[0.2cm]
c^{2,0}(g_1,g_2;A)&=&-Tr\lbrack \frac12(dg_2g_2^{-1})(g_1^{-1}dg_1)
(g_1^{-1}Ag_1)-\frac12(dg_2g_2^{-1})(g_1^{-1}Ag_1)
(g_1^{-1}dg_1)\rbrack ,\\[0.2cm]
c^{3,0}(g_1,g_2,g_3)&=&c^{2,0}(g_2,g_3,g_1^{-1}dg_1).
\end{eqnarray*}
In the above \(F_A\) is abbreviated to \(F\).

\vspace{0.2cm}

The {\it Chern-Simons form}  on \(N\) is by definition.   
\begin{equation}
c^{0,2}(A)=Tr(\,A F^2-\frac12 A^3F+\frac1{10} A^5\,),\qquad  A\in{\cal A}(N)\quad F=F_A.
\end{equation}  
The variation of the Chern-Simons form along the \({\cal G}(N)\)-orbit is given by (\ref{WZdescent2}):  
\begin{equation}
c^{0,2}( g\cdot  A)-c^{0,2}( A)=d\, c^{1,1}( g,\, A)+c^{1,2}( g), \quad g\in {\cal G}(N),
\label{eq:variation} 
\end{equation}

\subsection{ Fundamental formula}

Let  \(M\) be a connected compact four-manifold that is the boundary of an oriented 5-dimensional manifold \(N\); \(\partial N=M\). Let \(G\)  denote the Lie group \(SU(n)\) with \(n\geq 3\).    We have \(\pi_4(G)=1\).
We consider the trivial bundle \({\bf P}=N\times G\) and its restriction \(P\) to \(M\);\(\,P=M\times G\).   \({\cal A}(N)\) and \({\cal G}(N)=L^2_{s+\frac12}(N,G)\) denote the space of connections and the group of gauge transformations on \(N\) respectively.  
Let  \({\cal A}(M)\) be the space of irreducible connections on \(M\) and   \({\cal G}(M)=L^2_s(M,G)\) the group of gauge transformations.

It is well known that, for a  \(G\)-principal bundle \(P\) on a manifold \(X\) and a  closed subset \(K\) of \(X\),  any connection \(A\) of \(P\) on \(K\) has an extension to a connection of \(P\) on \(X\).

We shall prove a lemma which plays a fundamental role to construct a pre-quantum line bundle over the moduli space of connections.   Let \(g\in{\cal G}(M)\).    Since  \(\pi_4(G)=1\),  
 \(g\in {\cal G}(M)\) has an extension  
\({\bf g}\in {\cal G}(N)\), \({\bf g}\vert M=g\) and since \(H^5(G,{\bf Z})\simeq {\rm Z}\) the extension is unique modulo \({\rm Z}\).   We put 
\begin{equation}
C_5(g)=\frac{i}{24\pi^3}\int_{N}c^{1,2}({\bf g})= \frac{i}{240\pi^3}\int_{N}Tr(d{\bf g}\cdot {\bf g}^{-1})^5.
\end{equation}   
\(C_5(g)\) is well defined modulo \({\rm Z}\) independently of the
extension.  

For \(A\in{\cal A}(M)\) and \(g\in{\cal G}(M)\)  we put
\begin{equation}
\Gamma(g,A)=\frac{i}{24\pi^3}\int_{M}c^{1,1}(g;A)+C_5(g).\label{defgamma}
\end{equation}
It is defined modulo \({\rm Z}\).     

\begin{lem}\label{Gammacocycle}
For \(A\in{\cal A}(M)\) and \(g\in {\cal G}(M)\) it holds that
\begin{equation}
\Gamma(fg,A)=
\Gamma(g,\,f\cdot A)+\Gamma(f,A),\qquad \mod {\rm Z}
\end{equation}
\end{lem}
{\it proof.}~~
  Let  \({\bf A}\in{\cal A}(N)\) and \(
{\bf g}\in{\cal G}(N)\) be extensions of \(A\) and \(g\) respectively.   
Then from (\ref{eq:variation}) we have
\begin{equation}
(\delta c^{0,2})({\bf g},{\bf A})=\,dc^{1,1}({\bf g};{\bf A})+c^{1,2}({\bf g}).
\end{equation}
Integration over \(N\) and Stokes' theorem yield 
\[\int_N\,(\delta c^{0,2})({\bf g},{\bf A})=\int_M\,c^{1,1}(g;A)+\int_N\,c^{1,2}({\bf g})=-24\pi^3i\,\Gamma(g,A).\]
Hence \[0=(\delta \Gamma)(f,g,A)=\Gamma(g,f\cdot A)+\Gamma(f, A)-\Gamma(fg, A)
.\]
\hfill\qed

\begin{rem}

 If  \({\bf A}\in{\cal A}(N)\) which extends  \(A\in{\cal A}(M)\) is a reducible connection then there is a \({\bf g}\in {\cal G}(N)\) such that \({\bf g}\neq 1\) and \({\bf g}\cdot{\bf A}={\bf A}\).    For  \((A,\,g={\bf g}\vert M)\) we have  \(\Gamma(g,A)\in {\rm Z}\).   In fact, since \(A\in {\cal A}(M)\) is irreducible \(g={\bf g}\vert M=1\), hence the first term of (\ref{defgamma}) vanishes and \({\bf g}\in {\cal G}(N/M)\), here \(N/M\) is the space obtaind by shrinking \(M\) to one point.  By virtue of the fact \(H^5(G,{\rm Z})\simeq{\rm Z}\) we have \(C_5( g)\in {\rm Z}\,\) and 
 \(\Gamma(g,A)=0 \mod {\rm Z}\).
\end{rem}

\subsection{Line bundle over \({\cal B}(M)\)  for  \(M=\partial N\) }

 In this section \(M\) is a connected compact four-manifold that is  the boundary of an oriented 5-dimensional manifold \(N\); \(\partial N=M\).   Let \({\cal A}(M)\),  \({\cal G}(M)\) and \({\cal B}(M)={\cal A}(M)/{\cal G}(M)\) as in the previous section.      We consider the \(U(1)-\)valued function on \({\cal A}(M)\times {\cal G}(M)\):
\begin{equation}
\Theta(g,A)=\exp 2\pi i\Gamma(g,A).
\end{equation}
Lemma \ref{Gammacocycle} yield 
the cocycle condition:
\begin{equation}
\Theta(g,A)\Theta(h,g\cdot A)=\Theta(gh,A).\label{eq:cocycle}
\end{equation}
Therefore if we define the action of \({\cal G}(M)\) on \({\cal A}(M)\times {\rm C}\) 
by 
\[\left(g, \,(A,c)\,\right)\longrightarrow \left(g\cdot A,\,\Theta(g,A)c\right),\]
we have a complex line bundle:
\begin{equation}
{\cal L}(M)={\cal A}(M)\times {\rm C}/{\cal G}(M)\longrightarrow {\cal B}(M)={\cal A}(M)/{\cal G}(M).
\end{equation}
\(\Theta\) being \(U(1)\)-valued, \({\cal L}(M)\) can be endowed with a hermitian metric.    
The associated \(U(1)\)-principal bundle is  
\begin{equation}
{\cal P}(M)={\cal A}(M)\times U(1)/{\cal G}(M)\stackrel{\pi}{\longrightarrow}  {\cal B}(M).
\end{equation}
In the following we shall investigate the connection on \({\cal P}(M)\).     
   First we define a splitting of the tangent space \(T_{(A,z)}(\,{\cal A}(M)\times U(1)\,)\) to the vertical subspace \(V_{(A,z)}=\{(0,c)\in T_A{\cal A}(M)\oplus i{\rm R}\}\) and the horizontal subspace.    The horizontal subspace is given by  
\[Hor_{(A,z)}=\{(a,c)\in T_A{\cal A}(M)\oplus i{\rm R};\quad c=\theta_A(a)\,\},\]
where \(\theta\) is the 1-form defined in (\ref{theta}).   
   Then a horizontal distribution in the tangent space \(T{\cal P}(M)\) is defined as the image of  \(Hor\) by the projection: \(\sigma_{\ast}: T(\,{\cal A}(M)\times U(1)\,)\longrightarrow T{\cal P}(M)\) which is derived from the projection \(\sigma: {\cal A}(M)\times U(1)\longrightarrow {\cal P}(M)\).    
  So \(\left(a\,,c=\theta_A(a)\right)\) and \(\left (g^{-1}ag\,,\,c+\,\widetilde d\,(\, \frac{1}{2\pi i}\log \Theta(g,\cdot )\,)_Aa\,\right)\) for any \(g\in {\cal G}(M)\)  have to be projected to the same horizontal tangent vector at \({\cal P}(M)_{[A,c]}\).     
That is,  we must have 
\[g\cdot\theta=\theta+\,\widetilde d\,(\, \frac{1}{2\pi i}\log\Theta(g,\cdot)\,).\]
This is proved in the following
\begin{lem}\label{connection}
\begin{equation}
\widetilde d\,(\,\frac{1}{2\pi i}\log\Theta(g,\,A)\,)\,a\,=\, \widetilde d\,\Gamma(g,A)\,a=(\delta\,\theta_A(a) )(g).
\end{equation}
\end{lem}
{\it Proof.}~~
Since \(\widetilde d\, C_5(g)=0\),  we have 
\begin{equation}\label{eq:gamma}
\begin{split}
\widetilde d \,\Gamma(g,\,A\,)a&=
\frac{i}{24\pi^3}\int_MTr[-V(aF+Fa)\\[0.2cm]
&\quad +\frac12(A^2V+AVA+VA^2+V^2A+AV^2+VAV+V^3)\,a\,],
\end{split}
\end{equation}
where we put \(V=dg\cdot g^{-1}\) and we used the relation
\begin{equation}\label{eq:vanish}
\begin{split}
d\,Tr[(VAa-VaA)]&=Tr[V(Ad_Aa+d_AaA)-V(aF+Fa)\\[0.2cm]
&\qquad +(V^2A+2AVA+AV^2)a\,].
\end{split}
\end{equation}
On the other hand 
\begin{equation}
\begin{split}
\left(\delta\theta_A(a)\right)(g)&=g\cdot\theta_A(a)-\theta_A(a)\\[0.2cm]
&=-\frac{i}{24\pi^3}\int_M\,Tr[VFa+FVa -\frac12((A+V)^3a-A^3a)],
\end{split}
\end{equation}
which is equal to the right hand side of (\ref{eq:gamma}).
Therefore
\begin{equation}
(\delta\,\theta_A(a))(g)= \widetilde d\,\Gamma(g,A)(a).
\end{equation}
\hfill\qed 

From this lemma we see that the the connection form on \({\cal P}\) is given by
\begin{equation}
\widetilde d\phi\,-\theta\,,
\end{equation}
where \(\phi\) is the angle coordinate of \(U(1)\).\\
  
We have obtained the following theorem.
\begin{thm}~There exists a hermitian line bundle  with connection \({\cal L}(M)\longrightarrow {\cal B}(M)\).
\end{thm}

Let \(\kappa\) be the curvature form of \(\theta\).   The 2-form \(\omega=i\widetilde d\theta\) on \({\cal A}(M)\) does not descend to \({\cal B}(M)\) unless \(\omega\) is a flat connection, but restricted to the subspace \(H^0_A\,\subset T_A{\cal A}(M)\),  (1.10), it gives the value of \(\kappa\):
\begin{equation}
\kappa_{[A]}([a],[b])\,=\,-i\omega_A(a,b)=-i\omega^0_A(a,b),
\end{equation}
for \([a],\,[b]\in T_A{\cal B}(M)\) and \(a,b\,\in H^0_A\).

\subsection{The Chern-Simons functional on \(N\)}

 The manifold \(M=\partial N\) is the same as in the previous part.   
Let \({\cal P}(N)\) be the pullback of the \(U(1)\)-principal bundle \({\cal P}(M)\longrightarrow {\cal B}(M)\) by the boundary restriction map 
\(\,r:\,{\cal B}(N)={\cal A}(N)/{\cal G}(N)\longrightarrow {\cal B}(M)\,\).   
\begin{equation}
{\cal P}(N)=r^{\ast}{\cal P}(M)={\cal A}(N)\times_{r^{\ast}\Theta}U(1).
\end{equation}
Then the connection \(\theta\) pulls back to give a connection \(\theta_N\) on \({\cal P}(N)\).    
We shall investigate the section of \({\cal P}(N)\) induced by the Chern-Simons functional.    In the sequel we denote connections on \(N\) by bold face letters and those on \(M\) by roman letters.

   The {\it Chern-Simons functional on} \({\cal A}(N)\) is defined by 
\begin{equation}
{\rm CS}({\bf A})=\frac{i}{24\pi^3}\int_N\, c^{0,2}({\bf A}),\qquad{\bf A}\in {\cal A}(N) .
\end{equation}
The equation (\ref{eq:variation}) implies 
\begin{eqnarray}
{\rm CS}(\bf g\cdot {\bf A})&=&{\rm CS}({\bf A})+\Gamma(g, A),\nonumber \\[0.2cm]
\exp\,2\pi i\,{\rm CS}({\bf g}\cdot {\bf A})&=&\Theta(g,A)\,\exp\,2\pi i\,{\rm CS}({\bf A})\,,
\label{eq:5cocycle}
\end{eqnarray}
where \(A\in {\cal A}(M)\) and \(g\in {\cal G}(M)\) are respectively the restriction of \({\bf A}\) and \(\bf g\) to the boundary \(M=\partial N\).
Then the map 
\[ s\,:\,{\cal A}(N)\longrightarrow {\cal A}(N)\times U(1)\]
given by 
\begin{equation}
s({\bf A})=\left({\bf A},\,\exp\,2\pi i\,{\rm CS}({\bf A})\right)
\end{equation}
is a section of 
\({\cal P}(N)\), that gives a trivialization of \({\cal P}(N)\).

\begin{prop}\label{eq:gradient}
The gradient vector field of the section \(s\) over \({\cal B}(N)\) is given by 
\[\,\ast\frac{i}{8\pi^3}{\bf F}_{{\bf A}}^2\,.\]
\end{prop}
{\it Proof.}~~
Let \(A=r{\bf A}\in{\cal B}(M)\).   
Let \( {\bf a}\in T_{{\bf A}}{\cal B}(N)\).   We denote the image of \({\bf a}\) by the boundary restriction map \(r\) by  \(a=r_{\ast}{\bf a}\in T_{A}{\cal B}(M)\).    
We have
\begin{equation*}
\begin{split}
(\,\widetilde d\,{\rm CS}({\bf A})){\bf a}&=\frac{i}{24\pi^3}\int_N\,Tr [ \,{\bf F}_{{\bf A}}^2{\bf a}+({\bf F}_{{\bf A}}{\bf A}+{\bf A}{\bf F}_{{\bf A}}-\frac12{\bf A}^3)d_{{\bf A}}{\bf a}\\
&\qquad -\frac12({\bf A}^2{\bf F}_{{\bf A}}+{\bf A}{\bf F}_{{\bf A}}{\bf A}+{\bf F}_{{\bf A}}{\bf A}^2-{\bf A}^4){\bf a}\, ].
\end{split}
\end{equation*}
On the other hand 
\begin{equation*}
\begin{split}
(\theta_N)_{{\bf A}}({\bf a})&=\theta_{A}(a)=\frac{-i}{24\pi^3}\int_M\,Tr[( AF_ A+ F_A A-\frac12A^3) a\,]\\
&=\frac{i}{24\pi^3}\int_N\,Tr[\, -2{\bf F}_{\bf A}^2{\bf a}-\frac12({\bf A}^2{\bf F}_{\bf A}+{\bf F}_{\bf A}{\bf A}^2+{\bf A}{\bf F}_{\bf A}{\bf A}-{\bf A}^4){\bf a}\\
&\qquad \qquad +({\bf A}{\bf F}_{\bf A}+{\bf F}_{\bf A}{\bf A}-\frac12{\bf A}^3)d_{\bf A}{\bf a}\,], 
\end{split}
\end{equation*}
Therefore we have 
\begin{equation}
\begin{split}
(\nabla s)_{{\bf A}}{\bf a}&=(\frac{1}{2\pi i}\widetilde d\log s\,-\,\theta_N\,)({\bf a})=(\widetilde d\,{\rm CS}({\bf A})){\bf a}-(\theta_N)_{{\bf A}}({\bf a})\\
&=
\frac{i}{24\pi^3}\int_N\,Tr[\, 3{\bf F}_{{\bf A}}^2{\bf a}]=\,\langle \,{\bf a}\,,\,\ast\frac{i}{8\pi^3}{\bf F}_{{\bf A}}^2\,\rangle_{L^2(N)}.
\end{split}
\end{equation}
\hfill\qed

The proposition says that connections  \({\bf A}\in {\cal A}(N)\) such that \({\bf F}_{{\bf A}}^2=0\) are the critical points of the gradient vector field of Chern-Simons functional \({\rm CS}\), in particular the flat connections  are critical points.   
For a flat connection \({\bf A}\in {\cal A}^{\flat}(N)\), the Chern-Simons functional becomes
\begin{equation}
{\rm CS}({\bf A})=\frac{i}{240\pi^3}\int_N\,Tr\left({\bf A}^5\right).
\end{equation}

\subsection{Line bundle over \({\cal B}(M)\) for a manifold  \(M\) with boundary }

Let \(\widehat M\) be a compact four-manifold that is the boundary of a five-manifold \(N\).     Let \(M\) be a connected four-dimensional submanifold of \(\widehat M\) with smooth boundary \(\partial M\).
   We consider the trivial \(G\)-principal bundle \(P=M\times G\), which is the restriction of the trivial bundle \(\widehat P=\widehat M \times G\) to \(M\).   Let \({\cal B}(M)={\cal A}(M)/{\cal G}_0(M)\) be the moduli space of irreducible connections over \(M\) .     
 We shall construct a line bundle with connection over \({\cal B}(M)\) whose curvature form is \(-i\,\omega\).   First we note that any connection \(A\in{\cal A}(M)\) is the resstriction of a connection \(\widehat A\in {\cal A}(\widehat M)\), and that the extension \(\widehat A\in {\cal A}(\widehat M)\) is irreducible.   In fact compare the holonomy groups defined by \(A=\widehat A\vert M\) and \(\widehat A\) over a point \(u=(x,a)\in M\times G\subset  \widehat M\times G\).   
Also we extend any \(g\in {\cal G}_0(M)\) across the boundary \(\partial M\) by defining it to be the identity transformation on \((\widehat M\setminus \stackrel{\circ}{M})\times G \).   Then \(\widehat g=g\vee 1^{\prime}\in {\cal G}(\widehat M)\).   Here \(1^{\prime}\) is the 
constant function; \(\widehat M\setminus  \stackrel{\circ}{M} \ni x\longrightarrow 1'(x)=1\in G\).    

We put, for  \(A\in{\cal A}(M)\) and \(g\in {\cal G}_0(M)\),  
\begin{eqnarray}
\Gamma_M(g;A)&=&\frac{i}{24\pi^3}\int_{M}c^{1,1}(g, A)+C_5(\widehat g),.\\
\nonumber\\
\Theta_M(g;A)&=&\exp 2\pi i\,\Gamma_M(g;A).\label{ThetaM}
\end{eqnarray}
These are defined mod {\rm Z}.

\begin{lem}
For \(A\in{\cal A}(M)\) and \(g\in {\cal G}_0(M)\) we have
\begin{equation}
\Gamma_M(fg,A)=
\Gamma_M(g,\,f\cdot A)+\Gamma_M(f,A),\qquad \mod {\rm Z}.\label{Gammacocycle2}
\end{equation}
\end{lem}
{\it proof.}~~
The lemma follows from Lemma \ref{Gammacocycle}.
\begin{equation*}
\Gamma_M(g,A)=\int_{M}c^{1,1}(g, A)+C_5(\widehat g)=\int_{\widehat M}c^{1,1}(\widehat g, \widehat A)+C_5(\widehat g)=\Gamma(\widehat g,\widehat A).\end{equation*}
\hfill\qed

It follows from (\ref{Gammacocycle2}) that 
 \(\Theta_M(g,A)\) satisfies the cocycle condition:
\[\Theta_M(g,\,A)\,\Theta_M(h,\,\,g\cdot A)=\Theta_M(gh,\,A),\qquad \mbox{for \(g,h\in {\cal G}_0(M)\)}.\]
So if we define the action of \({\cal G}_0(M)\) on \({\cal A}(M)\times {\rm C}\)  by 
\[\left(g, (A,c)\right)\longrightarrow \left(g\cdot A\,,\,\Theta_M(g, A)\,c\right),\]
 we have a hermitian line bundle on \({\cal B}(M)\) with the transition function \(\Theta_M(g,A)\):
\begin{equation}
{\cal L}(M)={\cal A}(M)\times {\rm C}/{\cal G}_0(M) \longrightarrow{\cal B}(M).
\end{equation}

The connection form on the associated principal \(U(1)\)-bundle is given by the formula (\ref{theta}).    In fact Lemma \ref{connection} is also valid for \(M\) with boundary, so 
by the same calculation as in section 3.3 we have the connection  \(\theta\) on \({\cal L}(M)\),   and  the curvature \(\kappa\), (3.19).    We have proved the following
\begin{thm}~There exists a hermitian line bundle  with connection \({\cal L}(M)\longrightarrow {\cal B}(M)\).   
\end{thm}

\section{Pre-quantization of the moduli space of flat connections}

We follow the definition by Guillemin et al. \cite{GGK} of {\it pre-quantization} of a manifold endowed with a closed 2-form.   For a manifold \(X\) endowed with a closed 2-form \(\sigma\), we call a  {\it pre-quantization} of \((X,\sigma)\)  a hermitian line bundle \(({\bf L},\,<\,,\,>)\) over \(X\) equiped with a hermitian connection \(\nabla\) whose curvature is \(\sigma\), equivalently a  pre-quantization of \((X,\sigma)\) is a principal \(U(1)\) bundle \(\pi:P\longrightarrow X\) and a connection form \(\theta\) on \(P\) with curvature \(\sigma\).

Let \(M\) be as in the previous section a submanifold of \(\widehat M\)  with smooth boundary which may be empty.  We remember that \(\widehat M\) is the boundary of an oriented five-manifold.    
  We know from Theorem \ref{bsymplectic} that if \(\partial M\neq\emptyset\) the moduli space \({\cal M}^{\flat}(M)\) of flat connections has the pre-symplectic structure \(\omega^{\flat}\) .     If \( \partial M=\emptyset\), that is, if 
\(M=\widehat M\),  \({\cal M}^{\flat}(\widehat M)\) is a finite dimensional manifold ( possibly with singularity ).    \({\cal M}^{\flat}(\widehat M)\) is one-point if \(\widehat M\) is simply connected.     

In the previous section we constructed the hermitian line bundle 
 \[ \pi:\,{\cal L}(M)\longrightarrow {\cal B}(M),\] 
with the transition function given by \( \Theta_M(g;A)=\exp 2\pi i\,\Gamma_M(g;A)\) and a hermitian connection \(\theta\) on it with curvature \(\kappa\).   We note that restricted to flat connections it holds that \(\kappa=-i\omega^{\flat}\) if \(\partial M\neq\emptyset\) and \(\kappa=0\) if \(M=\widehat M\).   
We shall investigate the restriction  of the line bundle \({\cal L}(M)\) to \({\cal M}^{\flat}(M)\subset{\cal B}(M)\).   

First we study the case \(\partial M\neq\emptyset\).   
   The restriction of \({\cal L}(M)\) to \({\cal M}^{\flat}(M)\) is the line bundle 
\begin{equation}
{\cal L}^{\flat}(M)={\cal A}^{\flat}(M)\times {\rm C}/{\cal G}_0(M)\longrightarrow {\cal M}^{\flat}(M).\end{equation}
   The connection on \({\cal L}^{\flat}(M)\) is given by the formula  (\ref{theta}) .    It becomes 
\begin{equation}
\theta_A(a)=\frac{i}{48\pi^3}\int_M\,Tr[ A^3\, a \,].\label{flattheta}
\end{equation}
The curvature  \(\kappa\) is given by 
\[\kappa_{[A]}(a,b)=
-i\, \omega^{\flat}_{[A]}(a,b)= 
 -i\, \omega^{\prime}_A(a,b)\,.\]
When \(\partial M=\emptyset\), that is, when \(M=\widehat M\),
 \begin{equation}
{\cal L}^{\flat}(\widehat M)={\cal A}^{\flat}(\widehat M)\times {\rm C}/{\cal G}(\widehat M)\longrightarrow {\cal M}^{\flat}(\widehat M).\end{equation}
The connection on \({\cal L}^{\flat}(\widehat M)\) is given by the same formula as in (\ref{flattheta}).    It is a flat connection,(3.19), 
\[\kappa_{[A]}(a,b)=
\,-i\, \omega^{0}_A(a,b)=0\,.\]

We have obtained the following theorem.
\begin{thm}~ 
There exists a pre-quantization of the moduli space  \(\left({\cal M}^{\flat}(M),\,\omega^{\flat}\right) \), that is,
\begin{enumerate}
\item
There exists a hermitian line bundle 
\({\cal L}^{\flat}(\widehat M)
\longrightarrow {\cal M}^{\flat}(\widehat M)\) with a flat connection.
\item
For a proper submanifold \(M\) of \(\widehat M \) 
 there exists a hermitian line bundle  with connection \({\cal L}^{\flat}(M)\longrightarrow {\cal M}^{\flat}(M)\),   whose curvature is equal to the pre-symplectic form  
\(-i\,\omega^{\flat}\).
\end{enumerate}
\end{thm}
We call \({\cal L}^{\flat}(M)\) pre-quantum line bundle over \({\cal M}^{\flat}(M)\).

If \(\widehat M\) is simply connected then \( {\cal M}^{\flat}(\widehat M)\) is one-point  and 
\({\cal L}^{\flat}(\widehat M)\simeq {\rm C}\).

\section{ The action of  \(\Omega^3_0G\) on \({\cal M}^{\flat}(D)\) and its lift to 
\({\cal L}^{\flat}(D)\) }

\subsection{ Abelian extension of the group \(\Omega^3_0G\)}

Let \(G=SU(n)\) with \(n\geq 3\).   
Let \(\Omega^3G\) be the set of smooth mappings from
\(S^3\) to \(G\) that are based at some point.   \(\Omega^3G\) is not connected but is divided into connected components by the degree.    We put 
\begin{equation}\Omega^3_0G=\{g\in \Omega^3G;\, \deg\, g=0\}.\end{equation}
J. Mickelsson gave a Lie group extension of \(\Omega^3_0G\) by 
the abelian group \(Map({\cal A}_3,U(1))\), where \({\cal A}_3\) is the space of connections on \(S^3\)~\cite{Mic}.       
In the following we shall explain it after~\cite{K,Mi, Mic}.

The oriented 4-dimensional disc with boundary \(S^3\) is denoted by \(D\).    We write  \(DG=Map(D,G)\), the set of smooth mappings from \(D \) to \(G\) based at  \(p_0\in \partial D=S^3\).     The restriction to \(S^3\) of  a \(f\in DG\) has degree \(0\); \(\,f\vert S^3\in \Omega^3_0G\).    We let  \(D_0G=\{f\in DG;\,f\vert S^3=1\}\).    
We define the action of \(D_0G\) on  \(DG \times Map({\cal A}_3,U(1))\) by 
\begin{equation}
h\cdot\,(f\,,\lambda\,)=\,\left(\,f\,h\,,\,\lambda(\cdot) \Theta_{D}(h\,,\,f^{ -1}df\,)\,\right),
 \end{equation}
where
\((f,\,\lambda\,)\in DG\times Map({\cal A}_3,U(1))\).   From (\ref{ThetaM}) we have 
\begin{equation}
\Theta_{D}(h\,,\,f^{ -1}df\,)=\exp 2\pi i\left(\frac{i}{24\pi^3}\int_{D}c^{2,1}(f,h)+C_5(h\vee 1^{\prime})\right).\label{ftheta}
\end{equation}
   We consider the quotient space by this action;
\begin{equation}
\widehat{\Omega G}=DG\times Map({\cal A}_3,U(1))/D_0G.
\end{equation}
The equivalence class of \((f, \lambda)\) is denoted by \([ f, \lambda]\).   
The projection \(\pi: \widehat{\Omega G}\longrightarrow \Omega^3_0G\) is defined by \(\pi([f\,,\,\lambda])= f\vert S^3\).   Then \(\widehat{\Omega G}\) becomes a principal bundle over \(\Omega^3_0G\) with the structure group \(Map({\cal A}_3,U(1))\).   
The transition function is \(\chi(f,g)=\Theta_D(f^{-1}g, \,f^{-1}df )\) for \(f,\,g\in DG\) such that \(f\vert S^3=g\vert S^3\).   
Here the \(U(1)\) valued function \(\chi(\,f\,,g\,)\) is considered as a constant function in \(Map({\cal A}_3,U(1))\).   As we see from the example  (\ref{OmegaD}) the associated line bundle to \(\widehat{\Omega G}\) is nothing but the restriction of the line bundle \(\pi: {\cal L}^{\flat}(D)\longrightarrow {\cal M}^{\flat}(D)\) to the set of pure gauges \(\{[f^{-1}df]\,:\,f\in DG\}\subset  {\cal M}^{\flat}(D)\).

We shall define a group structure on \(\widehat{\Omega G}\) so that \(
\widehat{\Omega G}\) becomess a Lie group extension of \(\Omega^3_0G\).
  We define the multiplication of two elements in \(DG\times Map({\cal A}_3,U(1))\) by
\begin{equation}
(f,\,\lambda)\bullet (g,\,\mu)=\left(fg\,
,\,\lambda(\cdot)\,\mu_{f}(\cdot)\exp\,2\pi i\,\gamma_{D}(\,f,g;\,\cdot)\,\right),
\end{equation}
where \[
\mu_{f}(A)=\mu\left ( (f\vert S^3)^{-1}A(f\vert S^3)+(f\vert S^3)^{-1}d(f\vert S^3)\right ),\]
and \(\gamma_D(f,g\,; A)\) is the Mickelsson's 
2-cocycle  on \(D\) which is given by the following formula.
\begin{eqnarray}\gamma_D(f,g\,; A)&=&
\frac{i}{24\pi^3}\int_D(\delta c^{1,1})(\,f,\,g\,;\,A)
\nonumber\\[0.2cm]
&=& \frac{i}{24\pi^3}\int_{S^3}c^{2,0}(f,g\,;\,A)+\frac{i}{24\pi^3}\int_D\,c^{2,1}(f,g).
\end{eqnarray}  
Then \(DG\times Map({\cal A}_3,U(1))\) is endowed with a group structure.    
The associative law follows from the relation \(\delta\,\gamma_D=0\,\), that is valid from the property  \(\delta^2=0\).    By some calculations based on the properties (\ref{WZdescent1}) and (\ref{WZdescent2}) we can prove 
that the set of elements  \((h,\,\exp 2\pi i\,C_5(h\vee 1')\,)\) with \(h\in D_0G\) forms a normal subgroup of \(DG\times Map({\cal A}_3,U(1))\), and that 
\( [f,\lambda] = [g,\mu] \) if and only if there is a \(h\in D_0G\) such that \((g,\mu )=(f,\lambda)\bullet (h,\exp\,2\pi i\,C_5(h\vee1'))\).    Hence  the group structure descends to the quotient space   \(\widehat{\Omega G}\).          
   The group \(Map({\cal A}_3,U(1))\) being embedded as a normal subgroup of   \(\widehat{\Omega G}\),     \(\,\widehat{\Omega G}\)  becomes  a group extension of \(\Omega^3_0G\) by the abelian group  \(Map({\cal A}_3,U(1))\).

\subsection{The action of \(\widehat{\Omega G}\) on pre-quantum line bundles}

 The action of the group of total gauge transformations \({\cal G}(D)\) on \({\cal A}(D)\) yields the action on \({\cal B}(D)={\cal A}(D)/{\cal G}_0(D)\).   Since  \({\cal G}(D)=DG\), \({\cal G}_0(D)=D_0G\) and \(DG/D_0G\simeq \Omega^3_0G\), we have the action of \( \Omega^3_0G\) on \({\cal B}(D)\).   The action  of \( \Omega^3_0G\) on \({\cal B}(D)\) does not lift to the action on the line bundle \({\cal L}(D)\).   The abelian extension \(\widehat{\Omega G}\) is needed to have the lift on  \({\cal L}(D)\).    

 For \(A\in{\cal A}(D)\) and  \(f\in DG\), we put
\begin{equation}
\beta_D (f,A)=\frac{i}{24\pi^3}\int_D\,c^{1,1}(f,\,A).
\end{equation}
Notice that the following relations hold:
\begin{eqnarray}
 \delta \beta_D&=&\gamma_D,\\[0.2cm]
 \Gamma_D(f,A)&=&\beta_D(f,A)+C_5(f\vee 1'), \qquad\mbox{ for \(f\in D_0G\) .}
 \end{eqnarray}    
 Then we define the action of \(DG\times  Map({\cal A}_3,U(1)) \) on \( {\cal A}(D)\times{\rm C}\)  by
\begin{equation}
(\,f,\,\lambda\,)\bullet (A,c)=
\left(\,f\cdot A\,,\,c\,\lambda(A\vert S^3)\exp 2\pi i\,\beta_D(f\,,\,A)\,\right).
\end{equation}
It is a right action.   
By virtue of the relation 
\(\gamma_D=\delta \beta_D\) we can verify 
\[
(g,\mu)\bullet \left(\,(f,\lambda)\bullet (A,c)\right)=\left((f,\lambda)\bullet (g,\mu)\right)\,\bullet \,(A,c),\]
for \((f,\lambda),\,(g,\mu)\in DG\times Map({\cal A}_3,U(1))\) and \((A,c)\in {\cal A}(D)\times {\rm C}\).   
Hence the action is certainly well defined.   In particular we have, for \(h\in D_0G\),  
\begin{equation}
(\,h,\,\exp 2\pi i\,C_5(h\vee 1')\,) \bullet (A,c)\,=\,
\left(\,h\cdot A\,,\,c\,\Theta_D(h\,,\,A)\,\right).
\end{equation}
From the definition of \({\cal L}(D)\) and the fact that \([\,h,\exp 2\pi i\,C_5(h\vee 1')\,]\) for \(h\in D_0G\) 
gives the unit element of \(\widehat {\Omega G}\) we see that the above action descends to the action of 
 \(\widehat {\Omega G}\) on \({\cal L}(D)\).   Thus we have proved the following theorem.  
\begin{thm}
The line bundle \({\cal L}(D)\) carries an action of  \(\widehat {\Omega G}\)  that is equivariant with respect to 
the action of \(\Omega^3_0G\) on the base space \({\cal B}(D)\).
\end{thm}

Now we shall investigate the reduction of \({\cal L}(D)\) by the action of \(\widehat{\Omega G}\).    We remember that  the base space 
\({\cal B}(D)\) is a \({\cal G}/{\cal G}_0\)-principal bundle over \({\cal C}(D)={\cal A}(D)/{\cal G}(D)\), hence 
is a \(\Omega^3_0G\)-principal bundle over \({\cal C}(D)\).  
   The reduction of \({\cal L}(D)\) by the action of \(\widehat{\Omega G}\) becomes a line bundle on \({\cal C}(D)\), as will be described in the following.     \\
Let \[{\cal K}(D)=
{\cal A}(D)\times{\rm C}/\widehat{\Omega G},\]
and let \(\pi:\,{\cal K}(D)\longrightarrow {\cal C}(D) \)
be the  projection induced from \(\pi:\,{\cal L}(D)\longrightarrow {\cal B}(D)\).
Then \(\pi:\,{\cal K}(D)\longrightarrow {\cal C}(D)\) becomes a line bundle with the structure group \(Map({\cal A}_3, U(1))\).   This is the reduction of \({\cal L}(D)\) by the action of \(\widehat{\Omega G}\) .   

  We saw in section 2.2 that the action of \({\cal G}(D)\) on 
\({\cal M}^{\flat}(D)\) is infinitesimally symplectic,   Hence  the action of \(\Omega^3_0G\) on 
\({\cal M}^{\flat}(D)\) is also infinitesimally symplectic.    This action is lifted to the pre-quantum line bundle \({\cal L}^{\flat}(D)\longrightarrow {\cal M}^{\flat}(D)\) by the action of \(\widehat {\Omega G}\).   In fact it is the action of 
 \(\widehat {\Omega G}\) restricted to \({\cal L}^{\flat}(D)\) .
\begin{thm}
The line bundle \({\cal L}^{\flat}(D)\) carries an action of  \(\widehat {\Omega G}\)  that is equivariant with respect to 
the infinitesimally symplectic action of \(\Omega^3_0G\) on the base space \({\cal M}^{\flat}(D)\).   
The reduction of  \({\cal L}^{\flat}(D)\) by \(\widehat{\Omega G}\) becomes the complex line \({\rm C}\).
\end{thm}
In fact the reduction of \({\cal L}^{\flat}(D)\) by the action of \(\widehat{\Omega G}\) is nothing but 
 the restriction of the line bundle \({\cal K}(D)\) to \({\cal N}^{\flat}(D)\subset {\cal C}(D)\).   
 But \( {\cal N}^{\,\flat}(D)={\cal A}^{\flat}(D)/{\cal G}(D)\) being one point the restriction becomes
  \[{\cal K}(D)\vert_{{\cal N}^{\flat}(D)}\simeq {\rm C}.\]

\end{document}